\documentclass[12pt]{article}

\usepackage{amsmath,amssymb,mathrsfs,eucal}

\usepackage[T1]{fontenc}

\newtheorem{thm}{Theorem}

\def\mod{\text{mod }}

\begin{document}

\title{On a class of Lebesgue-Ljunggren-Nagell type equations} 

\author{Andrzej D\k{a}browski, Nursena G\"unhan and G\"okhan Soydan} 

\date{}

\maketitle

{\it Abstract}.  Given odd, coprime integers $a$, $b$ ($a>0$),  we 
consider the Diophantine equation $ax^2+b^{2l}=4y^n$,   $x, y\in\Bbb Z$, $l \in \Bbb N$, $n$ odd prime,  
$\gcd(x,y)=1$. 
We completely solve the above Diophantine equation 
for $a\in\{7,11,19,43,67,163\}$, and $b$ a power of an odd prime, 
under the conditions  
$2^{n-1}b^l\not\equiv \pm 1(\mod a)$ and $\gcd(n,b)=1$. 
For other square-free integers $a>3$ and $b$ a power of an odd prime, 
we prove that the above Diophantine equation has  no solutions for all    
integers $x$, $y$ with ($\gcd(x,y)=1$), $l\in\mathbb{N}$  and all odd primes $n>3$, satisfying 
$2^{n-1}b^l\not\equiv \pm 1(\mod a)$, $\gcd(n,b)=1$, and 
$\gcd(n,h(-a))=1$, where $h(-a)$ denotes the class number of the imaginary quadratic field 
$\mathbb Q(\sqrt{-a})$.

\bigskip


\bigskip


\section{Introduction}

The Diophantine equation $x^2+C=y^n$ ($x\geq 1$, $y\geq 1$, $n\geq 3$) has 
a rich history. Lebesgue proved that this equation has no solution when $C=1$, and Cohn 
solved the equation for several values of $1\leq C\leq 100$. The remaining values of $C$ in the 
above range were covered by Mignotte and de Weger, and finally by Bugeaud, Mignotte 
and Siksek.  Barros in his PhD thesis considered the range $-100\leq C\leq -1$. 
Also, several authors (Abu Muriefah, Arif, D\k{a}browski, Le, Luca, Pink, Soydan, 
Togb\'e, Ulas,...) became interested in the case where only the prime factors of $C$ are 
specified. Surveys of these and many others topics can be found in \cite{AB} and \cite{BP}. 
Some people studied the more general equation $ax^2+C=2^iy^n$, $a>0$ and 
$i\leq 2$.

Given odd, coprime integers $a$, $b$ ($a>0$),  we 
consider the Diophantine equation 
\begin{equation}\label{maineq} 
ax^2+b^{2l}=4y^n,   \quad  x, y\in \Bbb Z,\,\, l, n \in \Bbb N, \, n \, odd \, \, prime,  \gcd(x,y)=1. 
\end{equation} 
If $a \equiv 1 \, \mod 4$, then reducing modulo $4$ we trivially obtain that the 
 equation \eqref{maineq} has no solution.

It is known (due to Ljunggren \cite{Lj}) that the Diophantine equation $ax^2+1=4y^n$, $n\geq 3$, 
has no positive solution with $y>1$ such that $a \equiv 3 (\mod 4)$ and the class number of the 
quadratic field $\mathbb Q(\sqrt{-a})$ is not divisible by $n$. When $a=3$, then $3x^2+1=4y^n$ has the 
only positive solution $(x,y)=(1,1)$.

\bigskip

As our first result, we completely solve the equation \eqref{maineq} 
for $a\in\{7,11,19,$ $43,67,163\}$, under the conditions  
$2^{n-1}b^l\not\equiv \pm 1(\mod a)$ and $\gcd(n,b)=1$.

\begin{thm}\label{thm.1} 
Fix $p\in\{7,11,19,43,67,163\}$ and $b= \pm q^r$, with $q$ an odd prime different from $p$ 
and $r\geq 1$. 

$(i)$ The Diophantine equation 
\begin{equation}\label{sec.eq}
px^2+b^{2l}=4y^n,\, l\in\mathbb{N},\, \gcd(x,y)=1
\end{equation}
has no solutions $(p,x,y,b,l,n)$ with integers 
$x$, $y$ and  primes $n>3$, satisfying the conditions 
$2^{n-1}b^l\not\equiv \pm 1(\mod p)$ and $\gcd(n,b)=1$.  

$(ii)$ If $n=3$ and $p\not= 7$, then the equation \eqref{sec.eq} has no solutions $(p,x,y,b,l,3)$ 
satisfying the conditions $4b^l\not\equiv \pm 1(\mod p)$ and $\gcd(3,b)=1$.   

$(iii)$ If $n=3$ and $p=7$, then the equation \eqref{sec.eq} leads to $6$ infinite families of solutions,  
corresponding to  solutions of  Pell-type equations \eqref{Pell.1},   \eqref{Pell.2}, \eqref{Pell.3},  
\eqref{Pell.4}, and satisfying the conditions $4b^l\not\equiv \pm 1(\mod 7)$ and $\gcd(3,b)=1$. 

\end{thm}

{\bf Remarks.}  
(i)  The Diophantine equation \eqref{sec.eq} has many solutions (infinitely many ?) satisfying 
the conditions  $2^{n-1}b^l\equiv \pm 1(\mod p)$ and $\gcd(n,b)=1$.  Examples include 

\noindent 
$(p,x,y,b,l,n) \in\{ (7,\pm 1,$ $2,\pm 11,1,5),(11,\pm 1,3,\pm 31,1,5)$, 
$(7,\pm 7,2,\pm 13,1,7),$ $(19,\pm 1,5,\pm 559,1,7)$, 
$(11, \pm 253, 3, \pm 67, 1,11), (19, \pm 2531, 5, \pm 8579, 1, 11)$, \\ 
$(7,\pm 1,2,\pm 181,1,13),  (11, \pm 1801, 3, \pm 21929, 1,17 ), 
(7,\pm 457, 2, \pm 797, 1, 19),  \\ 
(7, \pm 967, 2, \pm 5197, 1, 23)\}$.

(ii)  If $b$ is divisible by at least two different odd primes, then the Diophantine equation \eqref{sec.eq} 
may have solutions satisfying the conditions $2^{n-1}b^l\not\equiv \pm 1(\mod p)$. 
Examples include 

\noindent 
$(p,x,y,b,l,n) \in 
\{ (7,103820535541,4,10341108537,1,37)$, \\ 
$(7,4865,46,1320267,1,7)$,  $(19,315003,49,909715,1,7)$, \\  
$(19,581072253,49,3037108805,1,11) \}$. 

(iii) Write the equation \eqref{sec.eq} as $px^2+b^{2l}=4y(y^{(n-1)/2})^2$ (compare \cite[p.116]{Lj}).
Now using $4y=u^2+pv^2$, taking $u=\pm 1$, and multiplying the equation by $p$, we arrive at the 
equation 

\begin{equation}\label{Ljunggren} 
X^2 - p(1+pv^2)Y^2 = -pb^{2l}. 
\end{equation}
If $b=\pm 1$, we obtain the equation (7') in \cite{Lj}. Ljunggren used an old result by Mahler 
to deduce that, if $p>3$, then  \eqref{Ljunggren} has no solution with $Y>1$ such that any prime 
divisor of $Y$ divides $p(1+pv^2)$ as well. 

(iv) Question: may we extend Ljunggren's idea to prove non-existence of solutions of our equation 
for some $b^l$ ? 


\bigskip

For a family of positive square-free integers $a$ with $h(-a)>1$ we can prove the following result 
(a variant of the results by Bugeaud \cite{Bu} and Arif and Al-Ali \cite{AA} in a case of the 
equation $ax^2+b^{2l+1}=4y^n$). Let $h(-a)$ denote the class number of the imaginary 
quadratic field $\mathbb Q(\sqrt{-a})$.

\begin{thm}\label{thm.2} 
Fix a positive square-free integer $a$, different from $3$, $7$, $11$, $19$, $43$, $67$, $163$, 
and  $b= \pm q^r$, with $q$ an odd prime not dividing $a$ and $r\geq 1$.  
Then the Diophantine equation \eqref{maineq} has no solutions $(a,x,y,b,l,n)$,     
with integers $x$, $y$ and primes $n>3$ satisfying the conditions 
$\gcd(n,h(-a))=1$,  $2^{n-1}b^l\not\equiv \pm 1(\mod a)$, and $\gcd(n,b)=1$. 
\end{thm}

{\bf Remarks.} (i) There are a lot of positive square-free integers $a$ with $rad(h(-a)) | 6$ 
(hypothetically, infinitely many):  $18$ values of $a$ with $h(-a)=2$, $54$ values of $a$ 
with $h(-a)=4$, $31$ values of  $a$ with $h(-a)=6$, etc. Here $rad(m)$ denotes the radical 
of a positive integer $m$, i.e. the product of all prime divisors of $m$. 

(ii) For fixed $a$ and $b$ we can (in some cases) use MAGMA \cite{Magma} to solve the Diophantine equation $ax^2+b^{2l}=4y^3$ (applying \texttt{SIntegralPoints} subroutine of MAGMA to associated families of 
elliptic curves). In a general case, one can try to prove a variant of Dahmen's result \cite{Da} saying that 
the above equation has no solution for a positive proportion of $l$'s, not divisible by $3$. 

(iii) The following variant of a result by Laradji, Mignotte and Tzanakis
 (see \cite[Theorem 2.3]{LMT}) follows 
immediately from out Theorem 2 (note that always $h(-p) < p$). 
Let $p$, $q$ be odd primes with $p\equiv 3 (\mod 8)$ and $p>3$. Then the Diophantine equation 
$px^2+q^{2l}=4y^p$ has no solution $(x,y,l)$ with positive integers $x$, $y$, $l$ satisfying 
$\gcd(x,y)=1$. 

(iv) Dieulefait and Urroz \cite{DU} used the method of Galois representations 
attached to $\mathbb Q$-curves  to solve the Diophantine equation $3x^2+y^4=z^n$. 
The authors suggest that their method can be applied to solve this type of equations with $3$ replaced 
by other values of $a$. We expect that their method can be extended to the case $ax^2+y^4=4z^n$ 
with small $a$ as well.

(v) We can solve the Diophantine $ax^2+b^{2l}=4y^n$ for relatively small values of $a>0$ (at least) 
 in positive integers $x,y,l,n$, $\gcd(x,y)=1$, $n\geq 7$ a prime dividing $l$, by using 
the Bennett-Skinner strategy \cite{BS}.  We treat some examples in Section 3.  Let us also mention that 
the smallest positive integer $a$ with $h(-a)=7$ is $71$, 
and one needs to consider newforms of weight $2$ and level $10082$.   

(vi) Pink \cite{P}  used estimates for linear forms in two logarithms in the complex and the $p$-adic case, 
to give an explicit bound for the number of solutions of the Diophantine equation  
$x^2+ (p_1^{\alpha_1}\cdots p_s^{\alpha_s})^2 =2y^n$  in terms of $s$ and 
$\max\{p_1, \cdots, p_s\}$.  We can prove analogous result 
 concerning the equations $px^2 +  (p_1^{\alpha_1}\cdots p_s^{\alpha_s})^2 = 4y^n$, 
with $p\in\{7,11,19,43,67,163\}$.

\section{Proofs of Theorems \ref{thm.1} and \ref{thm.2}}

{\it Proof of Theorem \ref{thm.1}}.  

Below in the proof, $b$ is a power of an odd prime $q \not= p$. 

As the class number of $\mathbb Q(\sqrt{-p})$ with $p\in\{7,11,19,43,67,163\}$ is $1$, 
we have the following factorization 

$$
{b^l+x\sqrt{-p}\over 2} \cdot {b^l-x\sqrt{-p}\over 2} = y^n. 
$$
Now we have 

$$ 
{b^l+x\sqrt{-p}\over 2} = \left({u+v\sqrt{-p}\over 2}\right)^n,  
$$ 
where $u$, $v$ are odd rational integers. Note that necessarily $\gcd(u,v)=1$. 
Equating real parts   we get 

$$
2^{n-1}b^l = u\sum_{r=0}^{(n-1)/2}\binom{n}{2r}u^{n-2r-1}(-p)^rv^{2r}. 
$$
As $u$ is odd, its possible values are among divisors of $b^l$. 
Here, we assume that $2^{n-1}b^l\not\equiv \pm 1 (\mod p)$.

(i) If $u=\pm 1$, then 
$2^{n-1}b^l = \sum_{r=0}^{(n-1)/2}\binom{n}{2r}(-p)^rv^{2r}$, and 
in particular $2^{n-1}b^l\equiv \pm 1 (\mod p)$, a contradiction. 

(ii)  If $u \not= \pm 1, \pm b^l$, then $q$ divides $pvn$. 
Since $\gcd(u,v) = \gcd(p,b) =1$, then $q$ divides $n$, a contradiction. 

(iii) Assume $u=\pm b^l$.  Put 
$\alpha = {v\sqrt{p}+b^li\over 2}$. Then $(\alpha + \overline{\alpha})^2 = v^2p$, 
$\alpha \overline{\alpha} = {1\over 4}(v^2p + b^{2l})$, and $\alpha/\overline{\alpha}$ 
is not a root of unity. Hence $(\alpha,\overline{\alpha})$ is a Lehmer pair. 
Note that ${\alpha^n - \overline{\alpha}^n \over \alpha - \overline{\alpha}}=\pm 1$. 
On the other hand, using \cite{BHV-M} we obtain that  
${\alpha^n - \overline{\alpha}^n \over \alpha - \overline{\alpha}}$ has primitive divisors 
for $n=11$ and all primes $n>13$, and hence our equation has no solution for $n=11$ and 
for primes $n>13$. Let us consider the cases $n\in\{3,5,7,13\}$ separately. 
Let us stress that the data in \cite{BHV-M} are given for equivalence classes of $n$-detective 
Lehmer pairs:  two Lehmer pairs $(\alpha_1,\beta_1)$, $(\alpha_2,\beta_2)$ are equivalent 
(we write $(\alpha_1,\beta_1) \sim(\alpha_2,\beta_2)$)  
if $\alpha_1/\alpha_2 = \beta_1/\beta_2 \in\{\pm 1,\pm i\}$.

\bigskip 

$\underline{n=3}$.  According to \cite{BHV-M}, we have two possibilities: 
(a) ${v\sqrt{p}+ui \over 2} \sim {\sqrt{1+\lambda} + \sqrt{1-3\lambda} \over 2}$, $\lambda\not=1$, 
or 
(b) ${v\sqrt{p}+ui \over 2} \sim {\sqrt{3^k+\lambda} + \sqrt{3^k-3\lambda} \over 2}$, $k>0$, $3\nmid \lambda$. 

In the case (a) we have four subcases: 
(i) $1+\lambda=v^2p$ and $1-3\lambda=-u^2$ or (ii) $1+\lambda=-v^2p$ and $1-3\lambda=u^2$ 
or (iii) $1+\lambda=-u^2$ and $1-3\lambda=v^2p$ or (iv) $1+\lambda=u^2$ and $1-3\lambda=-v^2p$ . 

In the subcase (i) we obtain a contradiction 
reducing the second equation modulo $3$. 

In the subcase (ii) we obtain relation $u^2=3pv^2+4$.  If $p\not =7$, then reducing 
this equation modulo $8$, we obtain $1\equiv 5(\mod 8)$, a contradiction.  Now the case $p=7$ 
leads to Pell-type equation 
\begin{equation}\label{Pell.1}
u^2-21v^2=4.
\end{equation}
Using the assumption $u$ and $v$ are odd for our equation, any solution to \eqref{Pell.1} is given by 
$$\dfrac{u_t+v_t\sqrt{21}}{2}=\left(\dfrac{u_0+v_0\sqrt{21}}{2}\right)^t,$$   
where $(u_0,v_0)=(5,1)$ is minimal solution and $3\nmid t.$ 
Thus an infinite family of solutions of equation  \eqref{maineq} is given by 
$$(x_t,y_t,b_t^l,n)=
\left(\dfrac{3u_t^2v_t-7v_t^3}{4},\dfrac{7v_t^2+u_t^2}{4},\dfrac{u_t^3-21u_tv_t^2}{4},3\right).$$ 
(see \cite[Proposition 6.3.16]{Coh} for the details about the equation \eqref{Pell.1}).\\

In the subcase (iii), note that 
$4+3u^2=v^2p$, and hence $7\equiv 3(\mod 8)$ if $p\not=7$, a contradiction. 
If $p=7$, then we need to consider the Diophantine equation 
\begin{equation}\label{Pell.2}
7v^2-3u^2=4.
\end{equation}
Such an equation has $3$ infinite families of solutions 
$(v,u) \in \{(s+3r,s+7r), (-s+3r,s-7r), (4s+18r,6s+28r)\}$, where $s^2-21r^2=1$. But since $u$ and $v$ are odd, 
one gets $2$ infinite families of solutions $(v,u) \in \{(s+3r,s+7r), (-s+3r,s-7r)\}$ for the equation \eqref{Pell.2}. 
Any solution to the equation $s^2-21r^2=1$ is given by
$$
s_t+r_t\sqrt{21}=(s_0+r_0\sqrt{21})^t 
$$
where $(s_0,r_0)=(55,12)$ is minimal solution. Thus 2 infinite families of solutions of equation  \eqref{maineq} 
are given by

\begin{eqnarray*}
x_t&=&-3r_ts_t^2+63r_t^3+21s_tr_t^2-s_t^3\\
y_t&=&2s_t^2+14s_tr_t+28r_t^2\\
b_t^l&=&-5s_t^3-63s_t^3r_t-231s_tr_t^2-245r_t^3
\end{eqnarray*}
or
\begin{eqnarray*}
x_t&=&-3s_t^2r_t+s_t^3-21s_tr_t^2+63r_t^3\\
y_t&=&2s_t^2-14s_tr_t+28r_t^2\\
b_t^l&=&-5s_t^3+63s_t^2r_t-231s_tr_t^2+245r_t^3
\end{eqnarray*}
with $n=3$. (see \cite[Theorems 4.5.1, 4.5.2]{AnAn} for details about the equation \eqref{Pell.2}).

In the subcase (iv) note that $4+v^2p=3u^2$, and hence $7\equiv 3(\mod 8)$ if $p\not=7$, 
a contradiction. If $p=7$, then reducing  $4+7v^2=3u^2$ modulo $7$ we obtain 
$\square = -\square$, a contradiction.  

In the case (b) we have four subcases: 
(i) $3^k+\lambda=v^2p$ and $3^k-3\lambda=-u^2$ or (ii) $3^k+\lambda=-v^2p$ and $3^k-3\lambda=u^2$ 
or (iii)  $3^k+\lambda=-u^2$ and $3^k-3\lambda=v^2p$ or (iv) $3^k+\lambda=u^2$ and $3^k-3\lambda=-v^2p$.  

In the subcase (i) note that $3v^2p=u^2+4\cdot 3^k$, hence necessarily $k=1$ (otherwise 
$3 | \gcd(u,v)$). Therefore $v^2p=3t^2+4$, where $u=3t$.  If $p\not =7$, then reducing 
this equation modulo $8$, we obtain $3\equiv 7(\mod 8)$, a contradiction. If $p=7$, then we 
need to consider the Diophantine equation
\begin{equation}\label{Pell.3}
7v^2-3t^2=4. 
\end{equation}
As in subcase (iii) of part (a)
 Such an equation has $3$ infinite families of solutions 
$(v,u) \in \{(s+3r,3(s+7r)), (-s+3r,3(s-7r)), (4s+18r,3(6s+28r))\}$, where $s^2-21r^2=1$. 
But since $u$ and $v$ are odd, one obtains $2$ infinite families of solutions $(v,u) \in \{(s+3r,3(s+7r)), (-s+3r,3(s-7r))$ 
for the equation \eqref{Pell.3}. Any solution to the equation $s^2-21r^2=1$ is given by
$$
s_t+r_t\sqrt{21}=(s_0+r_0\sqrt{21})^t
$$
where $(s_0,r_0)=(55,12)$ is minimal solution. Thus 2 infinite families of solutions of equation \eqref{maineq} 
are given by
\begin{eqnarray*}
x_t&=&567s_tr_t^2+99r_ts_t^2+5s_t^3+945r_t^3\\
y_t&=&4s_t^2+42s_tr_t+126r_t^2\\
b_t^l&=&-9s_t^3-63s_t^2r_t+189s_tr_t^2+1323r_t^3
\end{eqnarray*}
or
\begin{eqnarray*}
x_t&=&-567s_tr_t^2+99r_ts_t^2-5s_t^3+945r_t^3\\
y_t&=&4s_t^2-42s_tr_t+126r_t^2\\
b_t^l&=&-9s_t^3+63s_t^2r_t+189s_tr_t^2-1323r_t^3
\end{eqnarray*}
with $n=3$.

In the subcase (ii) note that $4\cdot 3^k+3pv^2=u^2$, hence necessarily $k=1$ (otherwise 
$3 | \gcd(u,v)$). Therefore $4=3t^2-pv^2$, where $u=3t$. If $p\not =7$, then reducing 
this equation modulo $8$, we obtain $4\equiv 0(\mod 8)$, a contradiction. 
Now reducing $4=3t^2-7v^2$ modulo $7$, we obtain $\square = -\square$, a contradiction again. 

In the subcase (iii) note that $3t^2p=u^2+4$, where $v=3t$. Now reducing modulo $3$, we obtain 
a contradiction. 

In the subcase (iv) note that $u^2=3pt^2+4$, where $v=3t$.  If $p\not =7$, then reducing 
this equation modulo $8$, we obtain $1\equiv 5(\mod 8)$, a contradiction.  Now the case $p=7$ 
leads to Pell-type equation
\begin{equation}\label{Pell.4}
u^2-21t^2=4.
\end{equation}
Since $u$ and $v$ are odd, one gets that $t$ is odd for the equation \eqref{Pell.4}. So, any solution 
to \eqref{Pell.4}   is given by 
$$\dfrac{u_m+t_m\sqrt{21}}{2}=\left(\dfrac{u_0+t_0\sqrt{21}}{2}\right)^m,$$ 
where $(u_0,t_0)=(5,1)$ is minimal solution and $3\nmid m.$ 
Thus an infinite family of solutions of equation  \eqref{maineq} is given by
$$(x_m,y_m,b_m^l,n)=
\left(\dfrac{3u_m^2v_m-7v_m^3}{4},\dfrac{7v_m^2+u_m^2}{4},\dfrac{u_m^3-21u_mv_m^2}{4},3\right),$$ 
with $v_m=3t_m$.
\bigskip 

$\underline{n=5}$.  According to \cite{BHV-M}, we have two possibilities to consider: 
(a)  ${v\sqrt{p}+ui \over 2} \sim {\sqrt{F_{k-2\epsilon}} + \sqrt{F_{k-2\epsilon}-4F_k} \over 2}$, 
$k\geq 3$,  $\epsilon=\pm 1$, 
or 
(b) ${v\sqrt{p}+ui \over 2} \sim {\sqrt{L_{k-2\epsilon}} + \sqrt{L_{k-2\epsilon}-4F_k} \over 2}$, 
$k\not=1$, $\epsilon=\pm 1$. Here $F_m$ and $L_m$ denote $m$-th Fibonacci and Lucas 
number respectively. 

In the case (a) we have four subcases: 
(i) $v^2p=F_{k-2\epsilon}-4F_k$ and $-u^2=F_{k-2\epsilon}$ or 
(ii) $-v^2p=F_{k-2\epsilon}-4F_k$ and $u^2=F_{k-2\epsilon}$ or 
(iii) $v^2p=F_{k-2\epsilon}$ and $-u^2=F_{k-2\epsilon}-4F_k$ or 
 (iv) $-v^2p=F_{k-2\epsilon}$ and $u^2=F_{k-2\epsilon}-4F_k$. 

In the subcase (i) we obtain $F_{k-2\epsilon}=-u^2<0$, a contradiction. 

In the subcase (ii), due to the fundamental work by Ljunggren \cite{Lj0} \cite{Lj2} 
(see also \cite[Section 2]{BMS2}) 
we can find all solutions to the equation  $u^2=F_{k-2\epsilon}$ ($k\geq 3$,  $\epsilon=\pm 1$).  
Ljunggren has  proved that the only squares in the Fibonacci sequence are 
$F_0=0$, $F_1=F_2=1$, and $F_{12}=144$.  

The case $k-2\epsilon=1$ gives $k=3$, $\epsilon=1$, $u^2=1$, hence using the first 
equation from (ii) we obtain $-v^2p=-7$, i.e. $p=7$, $v^2=1$. This case gives the 
solution $(p,x,y,b^l,n) = (7,\pm 1,2,\pm 11,5)$ which contradicts with $2^{n-1}b^l\not\equiv \pm (\mod p)$.

The case $k-2\epsilon=2$, gives $k=4$, $\epsilon=1$, $u^2=1$, hence using the first 
equation from (ii) we obtain $-v^2p=-11$, i.e. $p=11$, $v^2=1$. This case gives the solution 
$(p,x,y,b^l,n) = (11,\pm 1,3,\pm 31,5)$, which is impossible since $2^{n-1}b^l\equiv \pm 1$ $(\mod p)$. 

The case $k-2\epsilon=12$, gives $k=14$, $\epsilon=1$, $u^2=144$ or 
$k=10$, $\epsilon=-1$, $u^2=144$. The first possibility gives 
$-v^2p=-2^2 \cdot 11 \cdot 31$, i.e. $p=11\cdot 31$, $v^2=4$. 
The second possibility gives $-v^2p=-2^2\cdot 19$, i.e. $p=19$, $v^2=4$. 
This case gives no solution satisfying $\gcd(x,b)=1$ (note that both $u$ and $v$ 
are even).

In the subcase (iii), let us note that 
$u^2 = -(F_{k-2\epsilon}-4F_k) = F_k+F_{k+2\epsilon}$. 
Now a short look at the paper by Luca and Patel (see \cite[Theorem 1]{LP},  and their 
calculations in Section 5) shows that $k=4$, $\epsilon=-1$ is the only possibility. 
But then $F_6 = 8 = v^2p$, a contradiction. 

In the subcase (iv) we obtain  $F_{k-2\epsilon}=-v^2p<0$, a contradiction.

In the case (b) we have four subcases: 
(i) $v^2p=L_{k-2\epsilon}-4L_k$ and $-u^2=L_{k-2\epsilon}$ 
or (ii) $-v^2p=L_{k-2\epsilon}-4L_k$ and $u^2=L_{k-2\epsilon}$ 
or (iii) $v^2p=L_{k-2\epsilon}$ and $-u^2=L_{k-2\epsilon}-4L_k$ 
or (iv)  $-v^2p=L_{k-2\epsilon}$ and $u^2=L_{k-2\epsilon}-4L_k$. 

In the subcase (i) we obtain $L_{k-2\epsilon}=-u^2<0$, a contradiction. 

In the subcase (ii) we can find all solutions to the equation $u^2=L_{k-2\epsilon}$, 
($k\not= 1$, $\epsilon=\pm 1$).  By the work by Cohn \cite{Cohn} we know all solutions: 
$L_1=1^2$ and $L_3=2^2$.  

The case $k-2\epsilon=1$, gives $k=3$, $\epsilon=1$, $u^2=1$, hence using the first 
equation from (ii) we obtain $-v^2p=-13$, i.e. $p=13$, $v^2=1$. 
The case $k-2\epsilon=3$, gives $k=5$, $\epsilon=1$, $u^2=4$, hence using the first 
equation from (ii) we obtain $-v^2p=-40$, i.e. $p=10$, $v^2=4$. 
None of these two cases lead to solution of our Diophantine equation ($13$ is congruent to 
$1$ modulo $4$, while $10$ is even).

In the subcase (iii), let us note that $L_k+L_{k-2\epsilon} = 5F_{k-\epsilon}$. 
Therefore we need to determine all $k$ such that $5F_{k-\epsilon}$ is a square. 
Again, the paper by Bugeaud, Mignotte and Siksek \cite{BMS} shows that 
the only possibility is $5F_5 = 5^2$. But then $v^2p=L_4=7$ in case 
$\epsilon=1$ or $v^2p=L_6=18$ in case $\epsilon=-1$. In the first case we 
obtain $p=7$, but then $u^2=65$, a contradiction.  
The second case gives a contradiction by trivial observation.   

In the subcase (iv) we obtain $L_{k-2\epsilon}=-v^2p<0$, a contradiction.

\bigskip 

$\underline{n=7}$. According to \cite{BHV-M} we have six equivalence classes of $7$-defective 
Lehmer pairs. Two of them, $({1-\sqrt{-7}\over 2}, {1+\sqrt{-7}\over 2})$ and 
 $({1-\sqrt{-19}\over 2}, {1+\sqrt{-19}\over 2})$, 
come from our Lehmer pairs, giving  
$(p,x,y,b^l,n)\in \{(7,\pm 7,2,\pm 13,7), (19,\pm 1,5,$ $\pm 559,7)\}$, which are impossible since 
$2^{n-1}b^l\equiv \pm 1$ $(\mod p)$. 

\bigskip 

$\underline{n=13}$. The unique $13$-detective equivalence class 
$({1-\sqrt{-7}\over 2}, {1+\sqrt{-7}\over 2})$ leads to the solution 
$(p,x,y,b^l,n) = (7,\pm 1,2,\pm 181,13)$, which contradicts with $2^{n-1}b^l\not\equiv \pm 1 (\mod p)$.

\bigskip 
{\it Proof of Theorem \ref{thm.2}}. 

In this case, thanks to \cite[Lemma 1]{Bilu}, we can follow the same lines as in the proof of Theorem 1 
for $n>3$.

\section{Method via Galois representations and modular forms}

We will consider the Diophantine equation 
$ax^2+b^{2k}=4y^n$, for $a\in\{3, 7, 11, 15\}$  
in positive integers $x, y, k, n$,  $\gcd(x,y)=1$,  $n\geq 7$ a prime dividing $k$. 
We will apply the 
Bennett-Skinner strategy \cite{BS}, in particular we will use the results we need from \cite{BS}. 
We can compute systems of Hecke eigenvalues for conjugacy classes of newforms using MAGMA 
(or use Stein's Modular Forms Database provided the level is $\leq 7248$). 

\bigskip

{\bf Remarks.} (a) If $a \equiv 3 \, \mod 8$, then $y$ is necessarily odd:  if $y$ is even, then reducing modulo $8$ 
we obtain that the left hand side is congruent to $4$ modulo $8$, while the right hand side is congruent 
to $0$ modulo $8$, a contradiction.  

(b) If $a \equiv 7 \, \mod 8$, then $y$ is necessarily even.

\bigskip

(i) The Diophantine equation  $3x^2+b^{2k}=4y^n$ has no solution $(x,y,k,n)$, $xy\not=1$, $n\geq 7$ 
prime dividing $k$.   

We will consider a more general Diophantine equation $X^n+4Y^n=3Z^2$ ($n\geq 7$ a prime) 
and use \cite{BS}.  We are in case (iii) of \cite[p.26]{BS},  hence $\alpha\in\{1, 2\}$. 
From Lemma 3.2 it follows, that we need to consider the newforms of weight $2$ and levels $N\in\{36, 72\}$.

a) There is only one newform of weight $2$ and level $36$, 
corresponding to an elliptic curve $E$ of conductor $36$ with complex multiplication by 
$\mathbb Q(\sqrt{-3})$. Here we will apply \cite[Subsection 4.4]{BS}, to prove that $ab=\pm 1$. 
Assume (a.a) that $ab\not= \pm 1$. Then (using Prop. 4.6 (b)) holds, hence if $n=7$ or $13$, $n$ splits in 
$K=\mathbb Q (\sqrt{-3})$ and $E(K)$ is infinite for all elliptic curves of conductor $2n$. One checks that 
both primes $7$ and $13$ split in $K$.  Now using Cremona's online tables we check, that all elliptic 
curves of conductor $126=2\times 7 \times 3^2$ have rank zero, and all elliptic curves of conductor 
$234=2 \times 13 \times 3^2$, which are quadratic twists by $3$ of quadratic curves of conductor $26$, 
have rank zero too. 

b) There is only one newform of weight $2$ and level $72$, corresponding to isogeny class of elliptic 
curves of conductor $72$, with $j$-invariant $u/3^v$, with $v>0$ and $u$ some non-zero 
integer prime to $3$. To eliminate such an elliptic curve we use \cite[Prop. 4.4]{BS}.

\bigskip 

(ii) The Diophantine equation    $7x^2+b^{2k}=4y^n$  has no solution $(x,y,k,n)$,  $n\geq 7$ 
prime dividing $k$.   

We need to consider the newforms of weight $2$ and level $N=98$. 

There are two Galois conjugacy classes of forms of weight $2$ and level $98$.  
We will use numbering as in Stein's tables: we have $a_3(f_2)=\pm \sqrt{2}$ and 
 we can use \cite[Prop. 4.3]{BS} to eliminate $f_2$. On the other hand, the form 
$f_1$ corresponds to an elliptic curve of conductor $98$, with $j$-invariant 
 $u/7^v$, with $v>0$ and $u$ some non-zero integer prime to $7$. To eliminate such 
an elliptic curves we use \cite[Prop. 4.4]{BS}.

\bigskip 

(iii)   The Diophantine equation    $11x^2+b^{2k}=4y^n$  has no solution $(x,y,k,n)$,  $n\geq 7$, 
$n\not=11, 13$ prime dividing $k$. 

We need to consider the newforms of weight $2$ and levels $N\in\{484, 968\}$.  

a)   There are five Galois conjugacy classes of forms of weight $2$ and level $484$.  
We have $a_3(f_1)=1$,  $a_3(f_4)={1\pm \sqrt{33}\over 2}$, and we can use \cite[Prop. 4.3]{BS} 
to eliminate $f_1$ and $f_4$. To eliminate $f_2$ and $f_3$ we need 
to consider coefficients $a_3$ and $a_5$:  we have 
$a_3(f_2)=a_3(f_3)={-3\pm \sqrt{5}\over 2}$ and 
$a_5(f_2)=a_5(f_3)={-1\pm \sqrt{5}\over 2}$ (we cannot avoid $n=29$ when 
considering only $a_3$).   Finally, $a_7(f_5)=\pm 2\sqrt{3}$, and we can use \cite[Prop. 4.3]{BS} 
to eliminate $f_5$ when $n\geq 7$ and $n\not=13$.

b)  There are fourteen Galois conjugacy classes of forms of weight $2$ and level $968$.  
We have $a_3(f_1)=-3$, $a_5(f_2)=a_5(f_3)=3$, $a_3(f_4)=a_3(f_5)=1$, 
$a_3(f_6)=a_3(f_7)=\pm 2\sqrt{5}$, $a_3(f_{10})={1\pm\sqrt{17}\over 2}$, 
and we can easily use \cite[Prop. 4.3]{BS} to eliminate $f_1,...,f_7$, and $f_{10}$. 
Now considering $a_3$ and $a_{13}$ for both newforms $f_8$, $f_9$, and $a_3$ and $a_5$ 
for both newforms $f_{13}$, $f_{14}$, and using \cite[Prop. 4.3]{BS}, we can 
eliminate these four forms when $n\geq 7$ and $n\not= 11$.

\bigskip 

(iv)  The Diophantine equation    $15x^2+b^{2k}=4y^n$  has no solution $(x,y,k,n)$,  $n\geq 7$ 
prime dividing $k$.  

We need to consider the newforms of weight $2$ and level $N=450$.  

 There are seven Galois conjugacy classes of forms of weight $2$ and level $450$. 
We have $a_{11}(f_3)=a_{11}(f_7)=3$,  and we can use \cite[Prop. 4.3]{BS} to eliminate $f_3$ and $f_7$. 
The forms $f_2$ and $f_6$ correspond to elliptic curves of conductor $450$ (named $C$ and 
$A$ respectively in Cremona's tables), with $j$-invariants $u/3^v$, with $v>0$ and $u$ some 
non-zero integer prime to $3$;  the forms $f_1$ and $f_5$ correspond to elliptic curves of conductor $450$ 
(named $F$ and $E$ respectively in Cremona's tables), with $j$-invariants $u/5^v$, with $v>0$ and $u$ some 
non-zero integer prime to $5$;  the form $f_4$ corresponds to elliptic curve $G$ of conductor $450$, 
with $j$-invariant $u/(3^v5^w)$, with $v,w>0$ and $u$ some 
non-zero integer prime to $15$. To eliminate all these elliptic curves we use \cite[Prop. 4.4]{BS}

\section*{Acknowledgments}

We would like to thank Professor Mike Bennett for providing us with some references. 
Andrzej D\k{a}browski would like to thank Professor  G\"okhan Soydan for inviting him 
to the {\it Friendly Workshop on Diophantine Equations and Related Problems} 
(Bursa Uluda\u{g} University, Bursa, July 6 - 8, 2019), where this work started. This workshop was partially supported by T\"{U}B\.{I}TAK (the Scientific and Technological Research Council of Turkey) under Project No: 117F287.

\bigskip 

Andrzej D\k{a}browski, Institute of Mathematics, University of Szczecin,
70-451 Szczecin, Poland,  E-mail: dabrowskiandrzej7@gmail.com and
andrzej.dabrowski@usz.edu.pl 

Nursena G\"unhan,  Department of Mathematics, Bursa Uluda\u{g} University, 16059 Bursa, Turkey,  
E-mail:  nursenagunhan@uludag.edu.tr

G\"okhan Soydan, Department of Mathematics, Bursa Uluda\u{g} University, 16059 Bursa, Turkey,  
E-mail:  gsoydan@uludag.edu.tr

\end{document}